\documentclass{article}
\usepackage[utf8]{inputenc}
\usepackage[final]{pdfpages}

\usepackage{latexsym,amssymb,amsmath}
\usepackage{epsfig}
\usepackage[english]{babel}

\usepackage{amsmath}
\newtheorem{theorem}{Theorem}[section]
\newtheorem{lemma}[theorem]{Lemma}
\newtheorem{proposition}[theorem]{Proposition}
\newtheorem{corollary}[theorem]{Corollary}
\newtheorem{definition}[theorem]{Definition}

\def\PP{\mathbb{P}}\def\AA{\mathbb{A}}\def\RR{\mathbb{R}}
\def\CC{\mathbb{C}}\def\HH{\mathbb{H}}\def\OO{\mathbb{O}}
\def\ZZ{\mathbb{Z}}\def\QQ{\mathbb{Q}}

\def\cO{{\mathcal O}}

\def\fg{{\mathfrak g}}\def\fso{{\mathfrak so}}

\title{The Cayley Grassmannian}
\author{Laurent MANIVEL}
\date{February 2016}

\begin{document}

\maketitle

\begin{abstract} We study the projective variety $CG$ parametrizing four dimensional
subalgebras of the complex octonions, which we call the Cayley Grassmannian. We prove that 
it is a spherical $G_2$-variety 
with only three orbits that we describe explicitely. Its cohomology ring has a basis 
of Schubert type classes and we determine the intersection product completely. 
\end{abstract}

\section{Introduction}

There exists only four real normed algebras up to isomorphism: the algebras of 
real and complex numbers, the Hamilton algebra of quaternions and the Cayley 
algebra of octonions. The Cayley-Dickson
doubling process allows to construct these algebras
iteratively as a chain
$$\RR \subset \CC \subset \HH \subset \OO .$$
The Cayley algebra  $\OO$ contains lots of subalgebras isomorphic to the 
Hamilton algebra. In fact its automorphism group $G_2=Aut(\OO)$ acts transitively
on the set of four dimensional subalgebras, which has a natural structure of compact 
manifold isomorphic with the homogeneous space 
$$ G_2/SO_3\times SO_3.$$

The main goal of this note is to describe in some details what happens over the 
complex numbers, that is, when we consider the complexified algebra of octonions. 
The set of four dimensional subalgebras of the complexified algebra of octonions
is a closed subvariety of the complex Grassmannian $G(4,8)$. Since all these subalgebras 
contain the unit element, we can focus on their imaginary part, which is parametrized 
by a closed subvariety $CG$ of the Grassmannian $G(3,7)$. An important difference with 
the real case is that the action of $G_2$ is no longer transitive. We will prove the following statement.

\begin{theorem}
The variety $CG$ is a smooth irreducible variety of dimension eight. The action of
$G_2$ on $CG$ is spherical, in particular quasi-homogeneous, and has only three orbits. 
\end{theorem}

A maximal torus in $G_2$ has fifteen fixed points in $CG$. There are infinitely
many invariant curves joining these points but we can nevertheless determine the 
equivariant cohomology ring of $CG$ by localization to the fixed points. Then 
we can easily derive the usual cohomology ring, which is isomorphic to the Chow ring. 

Interestingly the Betti numbers of $CG$ are the same as those of the Grassmannian
$G(2,6)$, and the presentation of the Chow ring that we obtain is quite similar. 
In particular the Picard number of $CG$ is one. We thus get a minimal
smooth compactification of the eight dimensional affine space, whose boundary 
we describe in some detail. 

Our description of the Chow ring is the following. Denote by $\sigma_1$ the hyperplane section. There is a Schubert type class $\sigma_2$, of degree two, such that:

\begin{theorem}The rational cohomology ring of $CG$ is
$$H^*(CG, \QQ) = \QQ [\sigma_1 , \sigma_2 ]/ \langle \sigma_1^5 - 5\sigma_1^3\sigma_2 + 
6\sigma_1 \sigma_2^2 , 16\sigma_2^3 -27\sigma_1^2 \sigma_2^2 + 9\sigma_1^4 \sigma_2\rangle .$$
\end{theorem}

The Picard group of $CG$ is cyclic, and the degree of $CG$ 
with respect to the (very) ample generator of the Picard group is $182$. We compute the generators of its homogeneous 
coordinate ring  thanks to the structure of its orbit closures.  
 We also compute the Chern classes of $CG$, and deduce
that its projective dual is a hypersurface of degree $17$ in $\PP^{28}$. 

In a subsequent paper we will compute the quantum cohomology of $CG$, which turns 
out to be semisimple. We also plan to describe its derived category which, if the 
the Dubrovin conjecture is correct, should be generated by an exceptional collection.

More generally, it would be extremely interesting to extend the huge amount of
information that we have about the quantum cohomology and the derived category of 
homogeneous spaces (although the picture is not complete yet) to the more general
setting of quasi-homogeneous, or more specifically, spherical varieties. Beyond the 
intrisic interest and beauty of the Cayley Grassmannian, we consider the present 
study as a tiny piece of this more general program.

\medskip\noindent
{\it Acknowledgements.} This paper was written  somewhere between 
the CRM (Montreal University) and the CIRGET (UQAM), during the Fall of 2014. It is a pleasure to
thank these institutions for the excellent working conditions and the 
stimulating environments. Many thanks also to
Steven Lu, Steven Boyer and Laurent Habsieger for their hospitality.

\section{Four dimensional subalgebras of $\OO$}

\subsection{A reminder about $G_2$}

From now on  we work exclusively over the field over complex numbers. We still denote by $\OO$ the 
algebra of octonions with complex coefficients. 
We will need two different points of views over the smallest of the exceptional
complex simple Lie groups. 
 
\subsubsection{$G_2$ and the octonions} 
One can describe $G_2$ as the automorphism group of the Cayley algebra of octonions, $G_2=Aut(\OO)$
(we use \cite{baez} as a convenient reference about the octonions). 
In particular it preserves the unit element
$e_0=1$, the multiplicative norm $q$, and therefore the space of imaginary octonions 
$Im(\OO)=e_0^{\perp}$.

The multiplication table of the octonions  can be encoded in the oriented 
Fano plane, whose vertices are in correspondence with the vectors $e_1, \ldots, 
e_7$ of an orthonormal basis of the imaginary octonions. 

\begin{center}
\setlength{\unitlength}{4mm}
\begin{picture}(30,12)(-15,-.7)
\put(-5.5,0){ \resizebox{!}{1.5in}{\includegraphics{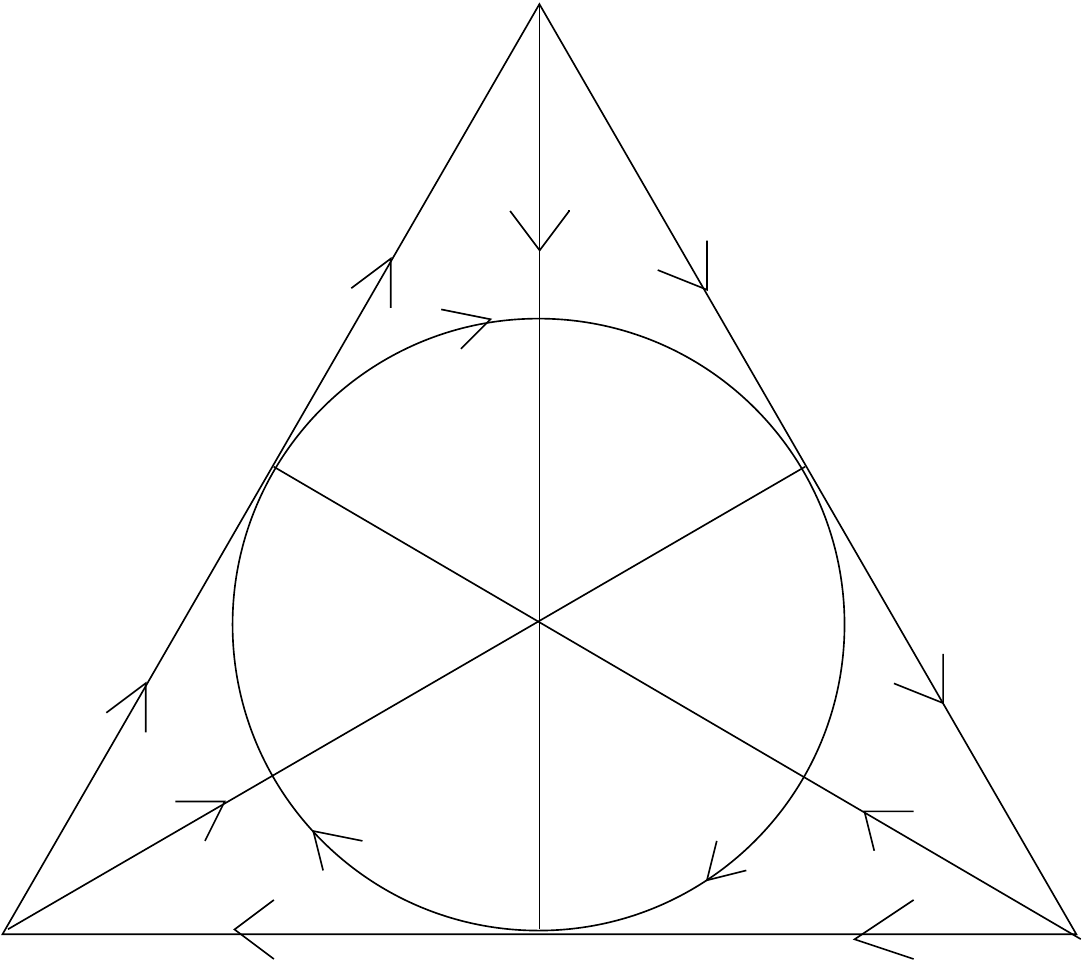}}}
\put(-3.8,4.7){$e_2$}\put(3.35,4.7){$e_6$}
\put(-.3,-.8){$e_4$}\put(-6.1,-.8){$e_1$}
\put(5.5,-.8){$e_5$}\put(-.2,9.9){$e_3$}\put(.9,3.2){$e_7$}
\end{picture} 
\end{center}

\centerline{\scriptsize{1. The multiplication table of the octonions}}

\medskip
This plane is made of seven lines (including the central circle) passing through 
seven points each. In order to multiply $e_i$ with $e_j$ for $i\ne j$, we need to find the line 
to which they both belong, and then $e_ie_j=\pm e_k$ for $e_k$ the third point
on this line; the sign is $+$ if the (cyclic) order $ijk$ of the vertices provides
the correct orientation of the line, and $-$ otherwise.  Finally $e_i^2=-1$ for any $i\ne 0$. 

The previous picture also defines a skew-symmetric three-form, given as the sum of 
the contributions of the seven (cyclically oriented) lines in the plane (note that 
the norm identifies $Im(\OO)$ with its dual):
$$\Omega = \sum_{\mathrm{lines}\;(ijk)}e_i\wedge e_j\wedge e_k.$$
In terms of the octonionic multiplication, this can be written as 
$$\Omega = \frac{1}{6}\sum_{ij}e_i\wedge e_j\wedge e_ie_j.$$

\subsubsection{The algebraic group $G_2$}
Recall that the algebraic group $G_2$ is a simply connected group of rank two and 
dimension fourteen. One of its fundamental representations is the adjoint representation
on the Lie algebra $\fg_2$, which is encoded in the root system also denoted $G_2$:

\setlength{\unitlength}{2mm}
\begin{picture}(20,37)(2,2)
\put(30,20){\vector(1,0){10}}
\put(30,20){\vector(1,2){4.5}}
\put(30,20){\vector(1,-2){4.5}}
\put(30,20){\vector(-1,0){10}}
\put(30,20){\vector(-1,2){4.5}}
\put(30,20){\vector(-1,-2){4.5}}

\thicklines
\put(30,20){\vector(0,1){17}}
\put(30,20){\vector(2,1){17}}
\put(30,20){\vector(-2,1){17}}
\put(30,20){\vector(0,-1){17}}
\put(30,20){\vector(2,-1){17}}
\put(30,20){\vector(-2,-1){17}}

\put(40.9,19.7){$\alpha_1$}
\put(10,28){$\alpha_2$}
\put(31,36){$\psi$}
\put(35.5,28){$\alpha$}\put(18,19.6){$\beta$}\put(35.5,11.5){$\gamma$}
\end{picture}

\centerline{\scriptsize{2. The  root system of type $G_2$}}

\medskip
This root system is made of two root systems of type $A_2$, formed by the short 
and the long roots respectively. The highest weight of the adjoint representation is
the highest root $\psi=3\alpha_1+2\alpha_2$.

The other fundamental representation is generated by the short roots. Its highest weight 
is the highest short root $\alpha=2\alpha_1+\alpha_2$; its weights are the six short roots 
and zero, all with multiplicity one, so that its dimension is seven; we denote this
representation by $V_7$. The action of $G_2$ on $V_7$ preserves a non degenerate 
quadratic form $q$, and also a skew-symmetric three form $\Omega$. (In fact we can 
define $G_2$ exactly that way: $\wedge^3V_7^*$ is a prehomogeneous vector
space under the action of $GL_7$,  and the stabilizer in $GL_7$ of a generic element 
$\Omega$ is a copy  of $G_2$.) The quadratic form $q$ and the skew-symmetric three form $\Omega$ can be related through the following relation (where the choice of the non zero constant $c$ is arbitrary):
$$q(x) \Theta
= c i(x)\Omega\wedge i(x)\Omega\wedge\Omega, $$
where  $i(x)$ is the contraction operator with $x\in V_7$, and $\Theta$ is a generator of
$\wedge^7V_7^*$ .

Consider the decomposition into weight spaces 
$$V_7=\CC u_0\oplus \CC u_{\alpha}\oplus \CC u_{-\alpha}\oplus 
\CC u_{\beta}\oplus \CC u_{-\beta}\oplus 
\CC u_{\gamma}\oplus \CC u_{-\gamma}.$$
Here we denoted the  six  short roots of $G_2$ as $\pm \alpha, \pm \beta,
\pm \gamma$, with $\alpha+\beta+\gamma=0$,  in order to
enhance the natural $S_3$-symmetry. The vector $u_0$ has non zero norm, 
the other weight vectors
are isotropic and we can normalize the octonionic quadratic form so that
its expression in the dual basis is 
$$q=v_0^2+v_{\alpha}v_{-\alpha}+v_{\beta}v_{-\beta}+v_{\gamma}v_{-\gamma}.$$ 
The  invariant three-form $\Omega$ up to scalar can then be normalized  as
$$\Omega=v_0\wedge v_{\alpha}\wedge v_{-\alpha}+v_0\wedge v_{\beta}\wedge v_{-\beta}+
v_0\wedge v_{\gamma}\wedge v_{-\gamma}+v_{\alpha}\wedge v_{\beta}\wedge v_{\gamma}+v_{-\alpha}\wedge v_{-\beta}\wedge v_{-\gamma}.$$
In order to reconcile the two approaches, note that $Im(\OO)=V_7$, the invariant 
norm being the restriction of the multiplicative norm. Moreover if in the previous 
expression of $\Omega$, we let $v_{\alpha}
=e_{\alpha}+e_{-\alpha}$, $v_{-\alpha}=e_{\alpha}-e_{-\alpha}$, and so on, 
it becomes a sum of seven terms that define a plane projective geometry. 
Moreover the octonionic product can be recovered directly from $\Omega$. 
Indeed, note that the induced map from $V_7\times V_7$ to $V_7$ sending 
$(x,y)$ to $Im(xy)$ is skew-symmetric; in terms of the three-form, 
$$Im(xy)=\iota(y\wedge x)\Omega.$$

A last remark is in order. Consider the two fundamental representations, and their 
projectivizations: each of them contains a unique closed $G_2$-orbit. In $\PP(V_7)$
this is simply the quadric $\QQ^5$ defined by the norm. In $\PP(\fg_2)$ this is by 
definition the {\it adjoint variety} of $G_2$, which we denote by $G_2^{ad}$ and 
is also five dimensional. Since $\fg_2\subset\fso_7=\wedge^2V_7$, the adjoint variety
is a subvariety of the Grassmannian $G(2,V_7)$. It was proved in \cite{lm1} that 
the adjoint variety parametrizes the {\it null-planes} inside $V_7=Im(\OO)$: those on 
which the octonionic product is identically zero. In particular, the following statement 
will be useful. 

\begin{proposition}
The action of $G_2=Aut(\OO)$ is transitive
\begin{enumerate}
\item on the set of non zero isotropic imaginary octonions; 
\item on the set of null-planes.
\end{enumerate}
\end{proposition}

\smallskip

\subsection{Non degenerate subalgebras}

Let $\AA=\CC 1\oplus Im(\AA)$ be a four dimensional subalgebra of the complexified
Cayley algebra.

\begin{definition}
We say that $\AA$ is non degenerate  when the restriction of the octonionic norm 
to $\AA$ is a non degenerate quadratic form. Otherwise we say that $\AA$ is degenerate.
\end{definition}

Suppose that $\AA$ is non degenerate.  Choose a non isotropic vector $e_1$ in 
$Im(\AA)$, and normalize it so that 
its norm is equal to one. Then $e_1^2=-e_1\bar{e_1}=-1$. Chose another norm one 
vector $e_2$ in $Im(\AA)$, orthogonal to $e_1$. Then $e_2^2=-1$. Let $e_3=e_1e_2
=-e_2e_1$. Recall that although the Cayley algebra is non associative, it is {\it alternative}, which
means that any subalgebra generated by 
two elements is associative. Hence  $e_1e_3=e_1(e_1e_2)=(e_1e_1)e_2=-e_2$
and similarly $e_2e_3=e_1$. 
Moreover, $1, e_1, e_2, e_3$ is a basis of $\AA$. Indeed, if this was not the case we would 
have a relation $e_3=x+ye_1+ze_2$; multiplying by $e_2$ on the left and on the right we 
would get $e_2=xe_1-y-ze_3=-xe_1+y-ze_3$, hence $x=y=0$;  multiplying by $e_1$ on the left and on the right
 we get similarly $x=z=0$; hence $e_3=0$, a contradiction since 
$e_3$ has norm one. 

This proves that  $\AA$ is isomorphic with the complexified Hamilton algebra. Moreover, 
since the norm is non degenerate on $\AA$ as well as on $\OO$, we can choose a norm
one vector $e_4$ orthogonal to $\AA$. Then letting $e_5=e_3e_4$, $e_6=e_2e_4$ and 
$e_7=-e_1e_4$, we can check that we get a basis $1, e_1, \ldots , e_7$ of $\OO$ whose 
multiplication table is the standard one. For example, the Moufang identities imply that 
$$ e_5e_6 = -(e_4e_3)(e_2e_4)=-e_4(e_3e_2)e_4=e_4e_1e_4=e_1,$$
and so on. We get the following statement:

\begin{lemma}
The automorphism group of $\OO$ acts transitively on the set of its four dimensional 
non degenerate subalgebras, which are all copies of the complexified Hamilton algebra.
\end{lemma}

Now we compute the stabilizer of a  non degenerate subalgebra of $\OO$, say $\HH$.

\begin{proposition}
The stabilizer of $\HH$ in $G_2=Aut(\OO)$ is the product $H=A_1^sA_1^l$ of the two copies 
of $SL_2$ in $G_2$ generated by its short roots and its long roots, respectively. 
\end{proposition}

{\it Proof.}
Let $H$ denote the stabilizer of $\HH$ in $G_2$.
It follows from the proof of the previous lemma that the map $r: H\rightarrow Aut(\HH)\simeq SL_2$ is surjective. Choosing a fixed $e_4$ we get a splitting $Aut(\HH)\hookrightarrow H$. Varying $e_4$ we see that the 
kernel $K$ of $r$ is in bijection with the space of unit vectors in $\HH^{\perp}$, yielding another copy of $SL_2$. 
The fact that
these two $SL_2$ correspond to $A_1^s$ and $A_1^l$ is straightforward. In fact we may
suppose that $Im(\HH)$ is generated by the three vectors $u_0,u_\alpha, u_{-\alpha}$.
Then the (short) $A_1$ generated by $\alpha$ preserves $\HH$ and acts non trivially
on it, while the (long) $A_1$ generated by the orthogonal roots acts trivially on $\HH$. 
$\Box$

\subsection{Degenerate subalgebras}

Now suppose that $\AA$ is degenerate. 

\begin{lemma}
$\AA$ contains a null plane $N$. 
\end{lemma}

{\it Proof.}  
Let $x$ be a non zero vector in the kernel of 
the quadratic form. This means that $x\in\AA\subset x^{\perp}$. In particular $x^2=-q(x)=0$.
Moreover, since $Re(xy)=q(x,\bar{y})$, the endomorphism of $\AA$ defined by left 
multiplication by $x$ stabilizes $Im(\AA)$. Its square is zero since $x^2=0$
and $\OO$ is alternative. Therefore there exists $y$ in $Im(\AA)$, independent of $x$,
such that $xy=0$, hence also $yx=0$ since $x$ and $y$, being orthogonal, anticommute. 
Note that this forces $y$ to be isotropic, since multiplying $xy=0$ on the right by 
$\bar{y}=-y$ we get $q(y)x=0$. Hence $y^2=0$ and $N=\langle x,y\rangle$ is a null-plane.  
$\Box$ 

\begin{lemma}
$Im(\AA)$ is contained in $N^{\perp}$. Conversely, any three-dimensional subspace
of  $N^{\perp}$ containing $N$ is the imaginary part of a degenerate subalgebra of $\OO$. 
\end{lemma}

{\it Proof.}  
Consider $z$ in $Im(\AA)-N$. Then the non zero vector $xz$ belongs to the kernel of left multiplication by $x$ in $\AA$.
Since this kernel is $N$,  there exist scalars $a$ and $b$ such that $xz=ax+by$.   Necessarily $b\ne 0$, 
otherwise we would get $(x-a1)z=0$ and then $z=0$ since $x-a1$ is not isotropic 
($a\ne 0$ since $z$ is not in $N$). Hence $y$ is a combination of $x$ and $xz$,
which easily implies that $y$ and $z$ anticommute. But this means they are orthogonal and we can conclude that  $z$ belongs to $N^{\perp}$. So $Im(\AA)\subset N^{\perp}$. 

In order to prove the second part of the statement we may use the transivity 
of $G_2$ on the set of null-planes. This allows us to suppose that $N=\langle x,y\rangle$
with $x=e_1+ie_2$ and $y=e_6+ie_7$. Then  $N^{\perp}=N\oplus 
\langle e_3,e_4,e_5\rangle$ and an explicit computation yields
$$\begin{array}{lll}
xe_3 = ix, & xe_4 = iy, & xe_5 = y, \\
ye_3=-iy,  & ye_4=ix,  & ye_5=-x.
\end{array}$$
This implies the claim. 
$\Box$ 
 
 \medskip
Consider as above $z$ in $Im(\AA)-N$. Since $z\in N^{\perp}$, we have only 
two possibilities. Either $z$ is not isotropic, so that the rank of the multiplicative 
norm on $Im(\AA)$ is equal to one; in this case $N$ is the only null-space contained
in $\AA$. Either $z$ is isotropic, so that the multiplicative norm is identically zero 
on $Im(\AA)$. But then the endomorphism of $Im(\AA)$ defined by left 
multiplication by $z$ has again rank one, and by the same argument as above 
$z$ is contained in another null-plane $N'$ of $Im(\AA)$. Let $u$ be a generator 
of $N\cap N'$. Then $u.Im(\AA)=0$. But the kernel of the left 
multiplication by $u$ on $\OO$ is precisely $u\OO$, which is a four dimensional
isotropic subspace of $\OO$. In particular it cannot be contained in $Im(\OO)$
and we conclude that necessarily
$$ Im(\AA)= u\OO\cap Im(\OO),$$
which is indeed the imaginary part of a subalgebra of $\OO$. 

\smallskip
We have finally classified the degenerate four dimensional  subalgebras of $\OO$.

\begin{proposition}
Let $\AA$ be a degenerate four dimensional  subalgebra of $\OO$. Then:
\begin{enumerate}
\item either $Im(\AA)$ is totally isotropic, there is a unique isotropic line 
$\ell$ such that $Im(\AA)= \ell\OO\cap Im(\OO),$ and each plane in $Im(\AA)$
containing $\ell$ is a null-plane,
\item  or $Im(\AA)$ is not totally isotropic, it contains a unique null-plane $N$ 
and it is contained in $N^{\perp}$.
\end{enumerate}
\end{proposition}
 
Up to the action of $G_2$, a four dimensional subalgebra of $\OO$ is thus 
isomorphic to one of the following:
$$\begin{array}{rcl}
H_0 & = & \langle e_0, e_1, e_2, e_3\rangle , \\
H_1 & = & \langle e_0, e_1+ie_2, e_6+ie_7, e_3\rangle , \\
H_2 & = & \langle e_0, e_1+ie_2, e_6+ie_7, e_4-ie_5\rangle .
\end{array}$$

\smallskip\noindent {\it Remark}. Of course over the real numbers, only the first type
is relevant: the multiplicative norm on the real Cayley algebra being positive 
definite,  it restricts on any four dimensional subalgebra to a positive 
definite multiplicative norm, so that any four dimensional subalgebra is a 
copy of the Hamilton algebra. 

\section{The Grassmannian embedding}

\subsection{Orbit structure} 

Let $CG$ be the set of four dimensional subalgebras of $\OO$. By considering only
the imaginary parts of these subalgebras we can consider $CG$
as a closed subvariety of the Grassmannian $G(3, V_7)$,  stable under the action 
of $G_2$ (it is closed because the limit of a family of subalgebras is certainly also
a subalgebra).  We can rephrase the previous discussion in the
following way. 

\begin{proposition}
The action of $G_2$ on $CG$ has only three orbits, consisting of:
\begin{enumerate}
\item the orbit $\cO_0\simeq G_2/A_1^sA_1^l$ of non degenerate subalgebras, of dimension eight, which is a spherical $G_2$-homogeneous space;
\item the orbit $\cO_1$ of degenerate but non isotropic subalgebras, which is 
fibered over the adjoint variety of $G_2$, the fibers being complements of smooth
conics inside projective planes; in particular its dimension is seven;
\item the orbit $\cO_2$ of isotropic subalgebras, which is isomorphic with the 
quadric $\QQ^5$ in $\PP(V_7)$.
\end{enumerate}
\end{proposition}

{\it Proof.} 
The only fact that we have not proved yet is that  $H=A_1^sA_1^l$ is a spherical 
subgroup of $G_2$. But this already appears in \cite{kr}.
$\Box$ 

\medskip
One can check directly that $\cO_0$ is dense in $CG$, or equivalently that 
$CG$ is irreducible. We will even show later on that $CG$ is smooth and irreducible.
 
We will denote by $H$ the orbit closure $\bar{\cO}_1=\cO_1\cup \cO_2$. 
It admits an equivariant resolution of singularities
$$\begin{array}{ccccc} 
 && \PP (N^{\perp}/N)=\tilde{H} && \\
 & \swarrow && \searrow & \\
 G_2^{ad} &&&&  H.
\end{array}$$
The multiplicative norm descends to the quotient bundle $N^{\perp}/N$,
defining a bundle of smooth conics over $G_2^{ad}$ whose total space 
is a hypersurface $E$ inside $\tilde{H}$ which is contracted to the 
closed orbit $\cO_2$ in $H$. In particular $H$ is singular along $\cO_2$,
with nodal singularities.  

\subsection{The three-form again}
We can also consider $CG$ as a subvariety of the dual
Grassmannian $G(4,7)$, by mapping a subalgebra $\AA$ of $\OO$ to its 
orthogonal complement $\AA^{\perp}$, which is contained in $Im(\OO)=V_7$. 

Let $T$ denote the tautological rank four vector bundle on $G(4,7)$. 
By the Borel-Weil theorem, the space of global sections of $\wedge^iT^*$ is
$$H^0(G(4,7),\wedge^iT^*)=\wedge^iV_7^* \quad \mathrm{for}\;\; 0\le i\le 4.$$ 

\begin{proposition}
Let $\Omega$ be a generic element of $\wedge^3V_7^*$. The zero locus 
of the corresponding global section of the rank four vector bundle 
$\wedge^3T^*$, is projectively isomorphic to $CG$.
\end{proposition}

\begin{corollary}
The variety $CG$ is a smooth Fano eightfold with canonical line bundle $K_{CG}=\cO_{CG}(-4)$.
\end{corollary}

{\it Proof.} 
Since the vector bundle $\wedge^3T^*$ is generated by global section, and since $\Omega$ is generic, 
the section it defines vanishes, if anywhere, on a smooth eight dimensional 
subvariety $Z$ of the Grassmannian, invariant under the $G_2$-action. By a direct
computation it vanishes at the three points defined by $H_0, H_1, H_2$, 
hence on their $G_2$-orbits, hence on the whole of $CG$. 
This implies that $CG$ is a 
connected component of $Z$. Finally, we check the connectedness of $Z$ with
the help of the Koszul complex
$$0\rightarrow \cO_G(-3)\rightarrow T(-2)\rightarrow \wedge^2T(-1)\rightarrow 
\wedge^3T\rightarrow \cO_G\rightarrow \cO_Z\rightarrow 0.$$
Indeed, a straightforward application of Bott's theorem yields $h^0(\cO_Z)=1$,
which is equivalent to the connectedness of $Z$.
$\Box$ 

\begin{corollary}\label{span}
The linear span of $CG$ is $\PP (S^2V_7)\subset \PP (\wedge^3V_7).$
\end{corollary}

Since $V_7$ is self-dual, $S^2V_7$ embeds equivariantly inside $End(V_7)$, which can be mapped to $\wedge^3V_7$  by sending an endomorphism $u$ to the three-form $\Omega(u.,.,.)$. Hence 
the embedding of $\PP (S^2V_7)$ inside $ \PP (\wedge^3V_7).$

Note that $S^2V_7$ decomposes as the direct sum of a one-dimensional trivial module and an
irreducible $G_2$-module. In the projectivization of the latter, the closed  $G_2$-orbit
is the quadric $\QQ_5$ (in its degree two Veronese re-embedding). This is in agreement 
with our description of the closed orbit $\cO_2$.

\begin{corollary}
$CG$ is projectively normal in $\PP (S^2V_7)$, with Hilbert polynomial
$$P_{CG}(k)= \frac{(k+1)(k+2)^2(k+3)}{2880}\left( 13(k+2)^4+7(k+2)^2+4\right) .$$
In particular the degree of $CG$ is $182$. \end{corollary}

{\it Proof.} 
Twisting the Koszul complex  by  $\cO_G(k)$ we get 
$$0\rightarrow \cO_G(k-3)\rightarrow T(k-2)\rightarrow \wedge^2T(k-1)\rightarrow 
\wedge^3T(k)\rightarrow \cO_G(k)\rightarrow \cO_{CG}(k)\rightarrow 0.$$
For $k\ge 0$ the first four vector bundles in this sequence have no higher cohomology,
hence $CG$ is projectively normal. For $k=1$, since $\wedge^3T(1)=T^*$ we get a sequence 
$$0\rightarrow H^0(G,T^*)=V_7^*\rightarrow H^0(G,\cO_G(1))=\wedge^3V_7\rightarrow  H^0(CG,\cO_{CG}(1))\rightarrow 0.$$
Since as $G_2$-modules $\wedge^3V_7=V_7^*\oplus S^2V_7$ (see \cite{lie}), 
we deduce Corollary \ref{span}.
Finally, the Borel-Weil theorem allows us to describe $H^0(CG,\cO_{CG}(k))$ for $k\ge 0$ as an alternate
sum of $GL_7$-modules whose dimensions are  given by the Weyl dimension formula, for example. 
Therefore we can deduce the Hilbert polynomial $P_{CG}(k)$, and then the degree of $CG$ by 
looking at its leading term.
$\Box$ 

\medskip
In particular there are $119$ quadrics vanihing on $CG$. We don't know whether they generate
the ideal of $CG$ or not. This would follow, after reduction to positive characteristics,  from the Frobenius splitting of a triple product of copies of $CG$, compatibly with its diagonals. Although $CG$ itself, being smooth and spherical, is certainly Frobenius split \cite{bi}, we do not know which powers of it inherite this property. 

\medskip
Since $CG$ is spherical, it is also multiplicity free, which means that the decomposition 
of its homogeneous coordinate ring decomposes into  $G_2$-modules without multiplicities. 
We can obtain this decomposition explicitely by using the divisor $H$ and the associated
exact sequence 
$$0\rightarrow \cO_{CG}(-1)\rightarrow \cO_{CG}\rightarrow \cO_H\rightarrow 0.$$
Twisting by  $\cO_{CG}(k)$ and taking global sections (recall that $CG$ is Fano with canonical
class $\cO_{CG}(-4)$, so that $\cO_{CG}(k-1)$ has no higher cohomology for any $k\ge 0$), 
we see that it is enough to determine the space of global sections of $\cO_H(k)$. For
this we can pull back to $\tilde{H}$ and make use of its projective bundle structure over
the adjoint variety. We deduce that 
$$H^0(H,\cO_H(k))=H^0(G_2^{ad}, S^k(\wedge^2N\wedge N^{\perp})^\vee).$$
In order to deduce the latter from the Borel-Weil theorem we need to understand
the bundle $\wedge^2N\wedge N^{\perp}=\det (N)\otimes N^{\perp}/N$ and 
the corresponding representation of the parabolic subgroup $P_{\alpha_2}$ of
$G_2$. For that, recall that  the 
adjoint variety $G_2^{ad}$ embeds into the Grassmannian $G(2,V_7)$. Moreover 
$V_7$ is generated by the short roots, in particular  the highest short root $\theta=2\alpha_1+\alpha_2$ coincides with the fundamental weight $\omega_1$, 
while the highest long root  $\psi=3\alpha_1+2\alpha_2$ is the fundamental weight
$\omega_2$. We deduce that the line of highest weight in $G_2^{ad}$ corresponds
to the plane $N=\langle e_{\alpha_1+\alpha_2}, e_{2\alpha_1+\alpha_2}\rangle$ in $V_7$. 
Then $N^{\perp}/N=\langle e_{-\alpha_1}, e_0, e_{\alpha_1}\rangle$, which 
corresponds to the three-dimensional irreducible representation of the $SL_2$ 
associated with $\alpha_1$. Taking its $k$-th power we obtain the direct sum of 
the irreducible representations of highest weights $(k-2j)\alpha_1$, for 
$j$ and $k-2j$ non negative. Since $\alpha_1=2\omega_1-\omega_2$, we 
conclude that $S^k(\wedge^2N\wedge N^{\perp})$ is the sum of the irreducible 
representations of highest weights $(k-2j)2\omega_1+2j\omega_2$. We finally get:

\begin{proposition}
The equivariant Hilbert series of $CG$ is
$$H_{CG}^{G_2}(t)=(1-tV_{\omega_0})^{-1}(1-tV_{2\omega_1})^{-1}(1-t^2V_{2\omega_2})^{-1},
$$
where $\omega_0$ denotes the trivial weight. Stated differently,
$$H^0(CG,\cO_{CG}(k))= \bigoplus_{i+2j\le k}V_{2i\omega_1+2j\omega_2}.$$
\end{proposition}

\section{Torus fixed points and localization}

Let $T$ be a maximal torus in $G_2$.

\subsection{Fixed points and Betti numbers}

\begin{proposition}
The variety $CG$ contains exactly $15$ fixed points of $T$.
\end{proposition}

{\it Proof.} 
The fixed points of $T$ in $CG$ are the intersection points of $CG$ with the projectivized 
weight spaces of $T$ in $\wedge^3V_7$. Up to signs and up to the symmetry in $\alpha, \beta, \gamma$
the weights are $2\alpha$ (six weights of multiplicity one), $\alpha$ (six weights of multiplicity three),
$\alpha-\beta$ (six weights of multiplicity one), and $0$ (with multiplicity five). 

A weight vector of weight $2\alpha$ is $u_\alpha\wedge u_\beta\wedge u_{-\gamma}$. The corresponding 
four dimensional space is $\langle u_0, u_\alpha, u_\beta, u_{-\gamma}\rangle$, on which the three-form
$\Omega$ vanishes identically. We thus get six $T$-fixed points of that type in $CG$.

A basis of weight vectors of weight $\alpha$ is $u_\alpha\wedge u_\beta\wedge u_{-\beta}$, 
$u_\alpha\wedge u_\gamma\wedge u_{-\gamma}$,  $u_0\wedge u_{-\beta}\wedge u_{-\gamma}$. No non
trivial combination of these vectors is completely decomposable, so the intersection of the 
projectivized weight space with $G=G(3,V_7)$ reduces to the three corresponding points. But none of
these is in $CG$, since they correspond respectively to the four dimensional spaces 
$\langle u_0, u_\alpha, u_\gamma, u_{-\gamma}\rangle$, $\langle u_0, u_\alpha, u_\beta, u_{-\beta}\rangle$
and $\langle u_\alpha, u_{-\alpha}, u_{-\beta}, u_{-\gamma}\rangle$, on which $\Omega$ does not vanish identically. 

A weight vector of weight $\alpha-\beta$ is $u_0\wedge u_\alpha\wedge u_{-\beta}$. The corresponding 
four dimensional space is $\langle u_\alpha, u_{-\beta}, u_\gamma, u_{-\gamma}\rangle$, on which the three-form
$\Omega$ vanishes identically. We thus get six new $T$-fixed points of that type in $CG$.

Finally, a basis of the zero weight space is $u_0\wedge u_{\alpha}\wedge u_{-\alpha}$, 
$u_0\wedge u_{\beta}\wedge u_{-\beta}$, $u_0\wedge u_{\gamma}\wedge u_{-\gamma}$, 
$u_{\alpha}\wedge u_{\beta}\wedge u_{\gamma}$, $u_{-\alpha}\wedge u_{-\beta}\wedge u_{-\gamma}$. 
Again it is easy to see that the intersection with $G$ consists in the five corresponding 
points only. The corresponding four dimensional spaces are $\langle u_\beta, u_{-\beta}, u_\gamma, u_{-\gamma}\rangle$
and the two others obtained by permuting  $\alpha, \beta, \gamma$, and 
$\langle u_0, u_{-\alpha}, u_{-\beta}, u_{-\gamma}\rangle$ and its opposite. 
Among these five only the first three are in $CG$. 
$\Box$ 

\medskip
In that situation we can use the Byalinicki-Birula decomposition theorem \cite{bb}: 
if we fix a generic
one-dimensional subtorus of $T$, it will again have the same fixed points as $T$ and the 
corresponding attracting sets will decompose $CG$ into fifteen strata, isomorphic to affine spaces. The cycle classes
of their closures then form a basis of the Chow ring,  or of the integer cohomology ring -- both 
are free and isomorphic to each other. By analogy with the case of the usual Grassmannians, we will call these cycle classes {\it Schubert classes}. 

In order to determine the Betti numbers of $CG$, there remains to compute 
the dimensions of the strata and count the numbers of strata of each dimension. This is done by listing the weights of the tangent spaces and counting the numbers of negative ones. We get:

\begin{proposition}
The Betti numbers of $CG$ are $1,1,2,2,3,2,2,1,1$.
\end{proposition}

For future use we will collect the weights of the tangent bundle to $CG$ at each of its fixed points. 
To obtain them we simply use the exact sequence 
$$0\rightarrow T_{CG}\rightarrow T_{G|CG}=T^*\otimes Q\rightarrow\wedge^3T^*\rightarrow 0.$$
The weights of $T$ and $Q$ are immediate to read at each fixed point. 

For simplicity we denote by $(x,y,z)$ the point $[u_x\wedge u_y\wedge u_z]$. The numbers on the 
leftmost column of the following table  refer to the labelling that we use in the sequel. \medskip

\begin{center}
\begin{tabular}{lll}
$0$\hspace*{1cm} & $(\alpha, \beta, -\gamma)$ \hspace*{1cm}  & $  \alpha, \beta, 2\alpha, 2\beta, \alpha-\gamma, \beta-\gamma, -\gamma, -\gamma$ \\
$5$ & $(\beta, \gamma, -\alpha)$   &   $ \alpha, \beta, 2\alpha, 2\beta, \alpha-\gamma, \beta-\gamma, -\gamma, -\gamma$ \\
$6'$ & $(\gamma, \alpha, -\beta)$   &  $  2\alpha, \alpha, \alpha-\beta, -\beta, -\beta, \gamma-\beta, \gamma, 2\gamma$ \\
$3$ & $(\alpha, -\beta, -\gamma)$  &   $-2\beta, -\beta, \alpha-\beta, \alpha, \alpha, \alpha-\gamma, -\gamma, -2\gamma$ \\
$2'$ & $(\beta, -\alpha, -\gamma) $ &   $ -2\alpha, -\alpha, \beta-\alpha, \beta, \beta, \beta-\gamma, -\gamma, -2\gamma$ \\
$8$ & $(\gamma, -\alpha, -\beta)$  &  $ -\alpha, -\beta, -2\alpha, -2\beta, \gamma-\alpha, \gamma-\beta, \gamma, \gamma$ \\ 
$5'$ & $(0,\alpha, -\beta) $  &  $ \alpha, \alpha-\beta, \alpha-\beta, -\beta, -\gamma, \alpha-\gamma, \gamma-\beta, \gamma$  \\
$2$ & $(0,\alpha, -\gamma)$   &   $ -\beta, \alpha-\beta, \alpha, \beta, \alpha-\gamma, \alpha-\gamma, \beta-\gamma, -\gamma$ \\ 
$3'$ & $(0,\beta, -\alpha) $  &  $ -\alpha, \beta-\alpha, \beta-\alpha, \beta, -\gamma, \beta-\gamma, \gamma-\alpha, \gamma$  \\
$1$ & $(0,\beta, -\gamma) $   &   $ -\alpha, \beta-\alpha, \alpha, \beta, \alpha-\gamma, \beta-\gamma, \beta-\gamma, -\gamma$ \\
$6$ & $(0,\gamma, -\alpha) $  &   $ \beta, \beta-\alpha, -\alpha, -\beta, \gamma-\alpha,  \gamma-\alpha, \gamma-\beta, \gamma$ \\
$7$ & $(0,\gamma, -\beta) $  &   $ \alpha, \alpha-\beta, -\alpha, -\beta, \gamma-\alpha,  \gamma-\beta, \gamma-\beta, \gamma$  \\  
$4''$ & $(0,\alpha, -\alpha)$    &  $  -\beta, \alpha-\beta, \beta-\alpha, \beta, \alpha-\gamma, \gamma-\alpha, \gamma, -\gamma$ \\
$4'$ & $(0,\beta, -\beta) $  &  $ -\alpha, \alpha, \alpha-\beta, \beta-\alpha, \gamma , \gamma-\beta, -\gamma, \beta-\gamma$ \\
$4$ & $(0,\gamma, -\gamma)$  &  $-\beta, -\alpha, \alpha, \beta, \gamma-\alpha, \gamma-\beta, \alpha-\gamma, \beta-\gamma$
  \end{tabular}
\end{center}

\medskip

Note that given our generic one dimensional torus, we can replace the attractive strata
by the repulsive ones. Their closures give two sets of varieties which meet transitively
at one or zero points when their dimensions add up to the dimension of $CG$, so that
their cohomology classes are Poincar\'e dual basis of the cohomology ring. We deduce 
that our Schubert basis is Poincar\'e self-dual up to a permutation.

\medskip 
Since there is a unique codimension one stratum, we can deduce:

\begin{corollary}
The Picard group of $CG$ is $Pic(CG)=\ZZ\cO_{CG}(1)$. 
\end{corollary}

The affine cell of maximal dimension has for 
complement the closure of the codimension one
stratum, which is not difficult to analyze. We leave to the 
reader the proof of the following statement. 

\begin{proposition}
Let $\ell\in\QQ^5$ define an isotropic line in $Im(\OO)$. 
Let 
$$X_1(\ell)=\{\AA\in CG, \; A\cap \ell\OO\ne 0\},$$
a special hyperplane section of $CG$.
Then the complement of $X_1(\ell)$ in $CG$ is an eight dimensional affine space.
\end{proposition}

Since, as we will prove below, the cohomology ring is generated by the Schubert classes of degrees one and two, we will complete 
the picture by describing the codimension two strata of the stratification, or rather their closures, that we denote by $X_2$ and $X'_2$. 

\begin{proposition}
There exists $\ell\in\QQ^5$ such that
$$X'_2=X'_2(\ell)=\{\AA\in CG, \; Im(\AA)\cap \ell\OO\ne 0\}.$$
There exists a null-plane $N$, containing $\ell$, such that 
$$X_2=X_2(N)=\{\AA\in CG, \; \dim (Im(\AA)\cap N^{\perp})\ge 2\}.$$
\end{proposition}

Note that these are both restrictions of special codimension two Schubert cycles on the usual Grassmannian.

\subsection{Equivariant cohomology classes}

Our Schubert classes define a basis of the cohomology ring over $\ZZ$, 
and also of the  equivariant cohomology ring over the polynomial ring 
$\Lambda=H_T^*(\mathrm{pt})=\CC[\alpha_1,\alpha_2]=\CC[\alpha,\beta,\gamma]/\langle \alpha+\beta+\gamma\rangle $. Recall (for example from
\cite{br}) that 
the equivariant cohomology classes can be defined in terms of the 
GKM-graph $\Gamma$ of the variety $CG$, whose vertices are the $T$-fixed 
points, and where there is an edge between two vertices when the 
corresponding $T$-fixed points can be joined by a $T$-equivariant 
curve. The graph $\Gamma$ is represented below, it has a nice 
$S_3\times\ZZ_2$-symmetry.

\bigskip\noindent {\it Remarks}. Although there are finitely many $T$-fixed 
points in $CG$, it is not true that there are only finitely many $T$-stable 
curves joining them. It is not difficult to find the exceptions to finiteness:
for each boundary triangle of the GKM graph $\Gamma$ there is a $\PP^2$ covered
by the $T$-stable curves joining the points corresponding to the vertices 
of the triangle; and for each of the three branches of the central star of
$\Gamma$ there is a $\PP^1\times\PP^1$ covered by  $T$-stable curves.

\medskip
The choice of a generic one-dimensional torus of $T$ breaks the symmetry of the picture,
each vertex being put in correspondence with a Schubert variety, obtained as the closure of the 
attractive set of the associated $T$-fixed point $p$. We will denote 
the Schubert  classes by $\sigma_i$, $\sigma'_i$, $\sigma''_i$ if necessary, where $i$
is the codimension. The correspondence with vertices in the GKM graph is encoded 
in the next figure, where for example   $\sigma'_i$ corresponds to the vertex labelled
$i'$, and so on. Note that Poincar\'e duality is given by the central symmetry of the 
graph, the three central vertices being fixed.

\vspace{1cm}

\setlength{\unitlength}{6mm}
\begin{picture}(20,19)(-2,-.3)
\put(0,4.33){\line(1,0){15}}
\put(0,13){\line(1,0){15}}
\put(0,4.33){\line(0,1){8.67}}
\put(15,4.33){\line(0,1){8.67}}
\qbezier(0,4.33)(3.75,2.16)(7.5,0)
\qbezier(0,13)(6.25,2.16)(7.5,0)
\qbezier(15,13)(8.75,2.16)(7.5,0)
\qbezier(15,4.33)(11.25,2.16)(7.5,0)
\qbezier(0,4.33)(5,13)(7.5,17.32)
\qbezier(0,13)(3.75,15.16)(7.5,17.32)
\qbezier(15,13)(11.25,15.16)(7.5,17.32)
\qbezier(15,4.33)(10,13)(7.5,17.32)
\qbezier(5,4.33)(6,10)(10,13)
\qbezier(5,13)(6,7.33)(10,4.33)
\qbezier(5,4.33)(9,7)(10,13)
\qbezier(5,13)(9,10.33)(10,4.33)
\qbezier(2.5,8.66)(7.5,11)(12.5,8.66)
\qbezier(2.5,8.66)(7.5,6.33)(12.5,8.66)

\put(-.15,4.13){$\bullet$}\put(4.85,4.13){$\bullet$}\put(9.85,4.13){$\bullet$}
\put(14.85,4.13){$\bullet$}\put(-.15,12.85){$\bullet$}\put(4.85,12.85){$\bullet$}
\put(9.85,12.85){$\bullet$}\put(14.85,12.85){$\bullet$}\put(7.35,-.13){$\bullet$}
\put(7.35,17.15){$\bullet$}\put(2.35,8.5){$\bullet$}\put(12.35,8.5){$\bullet$}
\put(6.2,8.5){$\bullet$}\put(7.9,9.6){$\bullet$}\put(7.9,7.35){$\bullet$}
\put(11.5,1.7){$\beta$}\put(3,1.8){$\alpha$}
\put(11.7,3.7){$\beta$}\put(2.6,3.8){$\alpha$}\put(6.8,3.7){$\small{\beta-\alpha}$}
\put(9.1,2.2){$\beta$}\put(5.4,2.2){$\alpha$}
\put(-.5,8.5){$\gamma$}\put(1,7){$\gamma$}\put(1,10){$\gamma$}
\put(15.2,8.5){$\gamma$}\put(13.6,7){$\gamma$}\put(13.6,10){$\gamma$}
\put(11.4,15.2){$\alpha$}\put(3,15){$\beta$}
\put(11.7,13.1){$\alpha$}\put(2.6,13.1){$\beta$}\put(6.8,13.1){$\small{\beta-\alpha}$}
\put(9.1,14.7){$\alpha$}\put(5.5,14.7){$\beta$}
\end{picture}

\centerline{\scriptsize{2. The GKM graph of $CG$}}

\vspace{1cm}

\medskip
Consider a Schubert variety $X$. By the classical localization theorems, its equivariant 
cohomology class can be represented by a map $f_X : \Gamma\rightarrow \Lambda$
with the following properties:
\begin{enumerate}
\item each polynomial $f_X(q)$ is a homogeneous polynomial, whose degree is equal to the codimension of $X$ in $CG$;
\item $f_X(q)=0$ if the $T$-fixed point corresponding to $q$ does not belong to $X$;
\item $f_X(p)$ is the product of the weights of the $T$-action on the normal space
to $X$ at $p$;
\item suppose that two vertices $q$ and $r$ are joined in $\Gamma$ by an edge, 
and let $y$ be the weight of the $T$-action on the tangent space at either one of 
the corresponding  $T$-fixed points. Then $f_X(q)-f_X(r)$ must be divisible by $y$.
\end{enumerate} 

One has more information a priori on the polynomials $f_X(q)$ but that turns out to be enough 
to determine them completely.  For example, if $X=CG$, we have $f_{CG}(q)=1$ for any
vertex $q$ of $\Gamma$. More interestingly if $X=H$, the unique codimension one
Schubert variety, then $f_H(q)=\omega(q)-\omega(0)$, where  $\omega(q)$ is the 
weight of the $T$-fixed point $q$ (which is a line in some weight space), and the 
vertex corresponding to the open strata is denoted by $0$. 

Now we use the following inductive strategy. Consider a Schubert 
variety $X$ of codimension $k$, corresponding to a vertex $p$ of $\Gamma$. 
Denote by  $Y_1,\ldots ,Y_m$ the codimension $k+1$ Schubert varieties 
contained in $X$, and by $q_1,\ldots ,q_m$ the corresponding vertices 
(note that $m$ is at most three, and most often two or one). By induction 
we know the functions $f_{Y_1},\ldots ,f_{Y_m}$. The function $f_X(f_H-f_H(p))$
is homogeneous of degree $k+1$, and is supported on the vertices corresponding to
Schubert varieties of codimension bigger than $k$. 

\vspace{1cm}

\setlength{\unitlength}{6mm}
\begin{picture}(19,18)(-2,-.3)
\put(0.3,4.33){\line(1,0){4.5}}\put(5.3,4.33){\line(1,0){4.4}}\put(10.3,4.33){\line(1,0){4.4}}
\put(0.3,13){\line(1,0){4.5}}\put(5.3,13){\line(1,0){4.4}}\put(10.3,13){\line(1,0){4.4}}
\put(0,4.66){\line(0,1){8.07}}
\put(15,4.66){\line(0,1){8.07}}
\qbezier(0.2,4.13)(3.75,2.16)(7.3,.2)
\qbezier(0.2,4.63)(1.25,6.46)(2.4,8.4)
\qbezier(2.7,9)(3.8,10.8)(4.8,12.7)
\qbezier(5.25,13.3)(6.3,15.15)(7.35,17.05)
\qbezier(0.2,12.8)(1.25,10.87)(2.3,9)
\qbezier(2.6,8.4)(3.8,6.53)(4.8,4.63)
\qbezier(5.15,4.1)(6.3,2.18)(7.35,.4)

\qbezier(14.8,4.13)(11.25,2.16)(7.7,.2)
\qbezier(14.8,4.63)(13.75,6.46)(12.6,8.4)
\qbezier(12.3,9)(11.2,10.8)(10.2,12.7)
\qbezier(9.75,13.3)(8.7,15.15)(7.65,17.05)
\qbezier(14.8,12.8)(13.75,10.87)(12.7,9.1)
\qbezier(12.4,8.4)(11.2,6.53)(10.2,4.63)
\qbezier(9.75,4.1)(8.7,2.18)(7.65,.4)

\qbezier(0.2,13.1)(3.75,15.16)(7.3,17.2)
\qbezier(14.8,13.1)(11.25,15.16)(7.7,17.2)

\qbezier(5,4.7)(5.3,7)(6.2,8.5)
\qbezier(6.5,9)(7.5,11)(9.8,12.9)
\qbezier(8.4,7.95)(9.5,10)(10,12.7)
\qbezier(5.25,4.4)(7,6)(7.9,7.3)
\qbezier(5,12.63)(5.3,10.33)(6.2,8.84)
\qbezier(5.25,12.9)(7,11.33)(7.9,10)
\qbezier(8.3,9.58)(9.5,7.33)(10,4.63)
\qbezier(6.5,8.33)(7.5,6.33)(9.8,4.44)
\qbezier(2.7,8.8)(5,9.7)(7.85,9.9)\qbezier(2.7,8.53)(5,7.63)(7.85,7.43)
\qbezier(8.3,9.9)(10.5,9.7)(12.3,8.8)\qbezier(8.3,7.43)(10.5,7.63)(12.3,8.53)

\put(-.15,4.13){$6'$}\put(4.85,4.13){$7$}\put(9.85,4.13){$6$}
\put(14.85,4.13){$5$}\put(-.15,12.85){$3$}\put(4.85,12.85){$2$}
\put(9.85,12.85){$1$}\put(14.85,12.85){$2'$}\put(7.35,-.13){$8$}
\put(7.35,17.15){$0$}\put(2.35,8.5){$5'$}\put(12.4,8.5){$3'$}
\put(6.2,8.5){$4$}\put(7.9,9.6){$4''$}\put(7.9,7.35){$4'$}
\end{picture}

\centerline{\scriptsize{3. Indexing Schubert classes}}

\vspace{1cm}

\bigskip

The equivariant cohomology 
class that it represents is then necessarily a combination of the classes of 
$Y_1,\ldots ,Y_m$ with constant coefficients. Otherwise said, there exists 
scalars $a_1,\ldots ,a_m$ such that 
$$f_X(q)(f_H(q)-f_H(p))=a_1f_{Y_1}(q)+\cdots +a_mf_{Y_m}(q)$$
for each $q$ im $\Gamma$. Letting $q=q_i$ we get 
$$f_X(q_i)(f_H(q_i)-f_H(p))=a_if_{Y_i}(q_i).$$
Since $f_H(q_i)-f_H(p)$ is never zero this determines $f_X(q_i)$ up to the constant 
$a_i$. But then the divisibility conditions imposed by the GKM 
graph allow to determine this constant (in most cases, but not always, it is enough 
for that to consider $f_X(q_i)-f_X(p)$). Once we have found $a_1,\ldots ,a_m$, 
the polynomial  $f_X(q)$ is completely determined for each $q$. 

\vspace{1cm}

\setlength{\unitlength}{6mm}
\begin{picture}(20,19)(-2,-.5)
\put(0.3,4.33){\line(1,0){2.65}}\put(5.3,4.33){\line(1,0){4.4}}\put(12,4.33){\line(1,0){2.6}}
\put(0.3,13){\line(1,0){3.5}}\put(6.3,13){\line(1,0){3.4}}\put(10.3,13){\line(1,0){4.4}}
\put(0,4.66){\line(0,1){8.07}}
\put(15,4.66){\line(0,1){8.07}}
\qbezier(0.2,4.13)(3.75,2.16)(7.3,.2)
\qbezier(0.2,4.63)(1.25,6.46)(2.4,8.4)
\qbezier(2.7,9)(3.8,10.8)(4.8,12.7)
\qbezier(5.3,13.4)(6.3,15.15)(7.35,17.05)
\qbezier(0.4,12.8)(1.25,10.87)(2.3,9)
\qbezier(2.6,8,4)(3.8,6.53)(4.8,4.63)
\qbezier(5.15,4.1)(6.3,2.18)(7.35,.4)

\qbezier(14.8,4.13)(11.25,2.16)(7.7,.2)
\qbezier(14.8,4.63)(13.75,6.46)(12.6,8.4)
\qbezier(12.3,9)(11.2,10.8)(10.2,12.7)
\qbezier(9.75,13.3)(8.7,15.15)(7.65,17.05)
\qbezier(14.8,12.8)(13.75,10.87)(12.7,9)
\qbezier(12.4,8,4)(11.2,6.53)(10.2,4.63)
\qbezier(9.75,4.1)(8.7,2.18)(7.65,.4)

\qbezier(0.4,13.1)(3.75,15.16)(7.3,17.2)
\qbezier(14.8,13.1)(11.25,15.16)(7.7,17.2)

\qbezier(5,4.7)(5.3,7)(6.2,8.5)
\qbezier(6.5,9)(7.5,11)(9.8,12.9)
\qbezier(8.3,7.75)(9.5,10)(10,12.7)
\qbezier(5.25,4.4)(7,6)(7.9,7.3)
\qbezier(5,12.63)(5.3,10.33)(6.2,8.84)
\qbezier(5.25,12.9)(7,11.33)(7.9,10)
\qbezier(8.3,9.58)(9.5,7.33)(10,4.63)
\qbezier(6.5,8.33)(7.5,6.33)(9.8,4.44)
\qbezier(2.8,8.8)(5,9.7)(7.85,9.9)\qbezier(2.8,8.53)(5,7.63)(7.85,7.43)
\qbezier(8.3,9.9)(10.5,9.7)(12.3,8.8)\qbezier(8.3,7.43)(10.5,7.63)(12.3,8.53)

\put(-2.15,4.13){$2\beta-2\gamma$}\put(3.1,4.13){$\beta-3\gamma$}\put(9.85,4.13){$\alpha-3\gamma$}
\put(14.85,4.13){$2\alpha-2\gamma$}\put(-.5,12.85){$2\beta$}\put(4.85,12.85){$\beta$}
\put(9.85,12.85){$\alpha$}\put(14.85,12.85){$2\alpha$}\put(6.6,-.33){$-4\gamma$}
\put(7.35,17.15){$0$}\put(.65,8.5){$2\beta-\gamma$}\put(12.35,8.5){$2\alpha-\gamma$}
\put(6.1,8.5){$-2\gamma$}\put(7.9,10.05){$-2\gamma$}\put(7.9,6.85){$-2\gamma$}

\end{picture}

\centerline{\scriptsize{4. The equivariant hyperplane class $\sigma_1$}}

\vspace{1cm}

\bigskip
Our procedure therefore effectively computes all the equivariant cohomology 
classes, by descending induction on the degree. For example, we have given below  the equivariant class $\sigma_2$,  which together 
with the hyperplane class generates the cohomology ring. 

\smallskip

Note that along the way we also compute the products of the equivariant Schubert classes
by the equivariant hyperplane class. Otherwise said, we get 
the equivariant Monk formula for $CG$ as a byproduct. 

\vspace{1cm}

\setlength{\unitlength}{6mm}
\begin{picture}(20,18)(-2,-.5)
\put(0.3,4.33){\line(1,0){2.65}}\put(6.1,4.33){\line(1,0){2.6}}\put(12,4.33){\line(1,0){2.6}}
\put(0.3,13){\line(1,0){3.6}}
\put(6.75,13){\line(1,0){2.95}}
\put(10.3,13){\line(1,0){4.4}}
\put(0,4.66){\line(0,1){8.07}}
\put(15,4.66){\line(0,1){8.07}}
\qbezier(0.2,4.13)(3.75,2.16)(7.3,.2)
\qbezier(0.2,4.63)(1.25,6.46)(2.4,8.4)
\qbezier(2.7,9.1)(3.8,10.8)(4.8,12.7)
\qbezier(5.3,13.4)(6.3,15.25)(7.35,17.15)
\qbezier(0.4,12.8)(1.25,10.87)(2.3,9)
\qbezier(2.7,8,3)(3.8,6.43)(4.8,4.53)
\qbezier(5.15,4.1)(6.3,2.18)(7.35,.4)

\qbezier(14.8,4.13)(11.25,2.16)(7.7,.2)
\qbezier(14.8,4.63)(13.75,6.46)(12.6,8.4)
\qbezier(12.3,9)(11.2,10.8)(10.2,12.7)
\qbezier(9.75,13.3)(8.7,15.15)(7.65,17.05)
\qbezier(14.8,12.8)(13.75,10.87)(12.7,9)
\qbezier(12.4,8,4)(11.2,6.53)(10.2,4.63)
\qbezier(9.75,4.1)(8.7,2.18)(7.65,.4)

\qbezier(0.4,13.1)(3.75,15.16)(7.3,17.2)
\qbezier(14.8,13.1)(11.25,15.16)(7.7,17.2)

\qbezier(5,4.7)(5.3,7)(6.2,8.5)
\qbezier(6.5,9)(7.5,11)(9.8,12.9)
\qbezier(8.3,7.75)(9.5,10)(10,12.7)
\qbezier(5.55,4.7)(7,6)(7.9,7.3)
\qbezier(5,12.63)(5.3,10.33)(6.2,8.84)
\qbezier(5.25,12.9)(7,11.33)(7.9,10)
\qbezier(8.3,9.58)(9.5,7.33)(10,4.63)
\qbezier(6.5,8.33)(7.5,6.33)(9.6,4.64)
\qbezier(2.9,9)(5,9.8)(7.85,9.9)\qbezier(2.9,8.33)(5,7.53)(7.85,7.43)
\qbezier(8.3,9.9)(10.5,9.7)(12.1,9)\qbezier(8.3,7.43)(10.5,7.63)(12.1,8.33)

\put(6,-.33){$4\gamma (\gamma-\beta)$}

\put(-3.05,4.13){$4\beta(\beta-\gamma)$}
\put(3.1,4.13){$2(\beta-\gamma)^2$}
\put(8.85,4.13){$\gamma(4\gamma-\beta)$}
\put(14.85,4.13){$2\gamma(\gamma-\alpha)$}

\put(.65,8.5){$\beta(4\beta-\gamma)$}
\put(4.9,8.5){$\beta(\beta-\gamma)$}
\put(11,8.5){$\gamma(\gamma-\alpha)$}

\put(7.9,6.85){$\gamma(4\gamma-\beta)$}
\put(7.9,10.05){$-3\beta\gamma$}

\put(-2.8,12.85){$2\beta(\beta-\alpha)$}
\put(4,12.85){$\beta(\beta-\alpha)$}
\put(9.85,12.85){$0$}
\put(14.85,12.85){$0$}

\put(7.35,17.15){$0$}

\end{picture}

\centerline{\scriptsize{4. The equivariant class $\sigma_2$}}

\vspace{1cm}

\subsection{The cohomology ring}

Modding out by 
non constant homogeneous polynomials we deduce the usual Monk formula, 
giving the product of a Schubert class by the hyperplane class. We synthetize
this formula in the following graph, which we call the Bruhat graph:

\setlength{\unitlength}{6mm}
\thicklines
\begin{picture}(20,10.3)(-1,-2.5)
\put(0,2){$\bullet$}\put(2,2){$\bullet$}
\put(8,2){$\bullet$}\put(14,2){$\bullet$}\put(16,2){$\bullet$}

\put(4,4){$\bullet$}\put(4,0){$\bullet$}\put(6,4){$\bullet$}\put(6,0){$\bullet$}
\put(10,4){$\bullet$}\put(10,0){$\bullet$}\put(12,4){$\bullet$}\put(12,0){$\bullet$}

\put(8,-2){$\bullet$}\put(8,6){$\bullet$}

\put(0.2,2.2){\line(1,0){2}}\put(4.2,4.1){\line(1,0){2}}\put(4.2,4.3){\line(1,0){2}}

\put(4.2,.2){\line(1,0){2}}
\put(10.2,4.1){\line(1,0){2}}\put(10.2,4.3){\line(1,0){2}}
\put(10.2,.2){\line(1,0){2}}\put(14.2,2.2){\line(1,0){2}}

\put(2.2,2.2){\line(1,1){2}}\put(2.2,2.2){\line(1,-1){2}}
\put(6.2,4.2){\line(1,1){2}}\put(6.2,4.2){\line(1,-1){2}}
\put(6.2,.1){\line(1,1){2}}\put(6.2,.1){\line(1,-1){2}}
\put(6.2,.3){\line(1,1){2}}\put(6.2,.3){\line(1,-1){2}}
\put(8.2,2.1){\line(1,-1){2}}\put(8.2,2.3){\line(1,-1){2}}
\put(8.2,6.2){\line(1,-1){2}}
\put(8.2,2.2){\line(1,1){2}}
\put(8.2,-1.7){\line(1,1){2}}\put(8.2,-1.9){\line(1,1){2}}
\put(12.2,4.2){\line(1,-1){2}}\put(12.2,.2){\line(1,1){2}}

\put(4.1,4.2){\line(1,-2){2}}\put(4.3,4.2){\line(1,-2){2}}
\put(4.1,.2){\line(1,2){2}}\put(4.3,.2){\line(1,2){2}}\put(4.2,.2){\line(1,2){2}}
\put(10.1,4.2){\line(1,-2){2}}\put(10.3,4.2){\line(1,-2){2}}\put(10.2,4.2){\line(1,-2){2}}
\put(10.1,.2){\line(1,2){2}}\put(10.3,.2){\line(1,2){2}}

\put(-.5,1.5){$\sigma_0$}\put(1.5,1.5){$\sigma_1$}\put(3.5,-.5){$\sigma_2$}
\put(3.5,4.6){$\sigma'_2$}\put(5.5,-.5){$\sigma_3$}\put(5.5,4.6){$\sigma'_3$}
\put(7.8,-2.5){$\sigma_4$}\put(7,2.1){$\sigma'_4$}\put(7.8,6.5){$\sigma''_4$}
\put(10.1,-.5){$\sigma_5$}\put(10.1,4.6){$\sigma'_5$}\put(12,-.5){$\sigma_6$}
\put(12,4.6){$\sigma'_6$}\put(14,1.5){$\sigma_7$}\put(16,1.5){$\sigma_8$}

\end{picture}

\centerline{\scriptsize{5. The Bruhat graph of $CG$}}

\bigskip 

\medskip
  This must be read as follows: for each Schubert class $\sigma$, the product 
of $\sigma$ by the hyperplane class is the sum of the Schubert classes connected to
it on the column immediately left, with coefficients equal to the number of edges
that connect them. (For example $\sigma_2H=\sigma_3+3\sigma'_3$ and 
$\sigma'_2H=2\sigma_3+2\sigma'_3$.)  We can deduce the degrees of all the Schubert classes, 
that we indicate on the following version of the Bruhat graph. Note that we recover the degree of $CG$ as $182=5^2+11^2+6^2.$

\setlength{\unitlength}{6mm}
\thicklines
\begin{picture}(20,10.3)(-1,-2.5)
\put(0,2){$\bullet$}\put(2,2){$\bullet$}
\put(8,2){$\bullet$}\put(14,2){$\bullet$}\put(16,2){$\bullet$}

\put(4,4){$\bullet$}\put(4,0){$\bullet$}\put(6,4){$\bullet$}\put(6,0){$\bullet$}
\put(10,4){$\bullet$}\put(10,0){$\bullet$}\put(12,4){$\bullet$}\put(12,0){$\bullet$}

\put(8,-2){$\bullet$}\put(8,6){$\bullet$}

\put(0.2,2.2){\line(1,0){2}}\put(4.2,4.1){\line(1,0){2}}\put(4.2,4.3){\line(1,0){2}}

\put(4.2,.2){\line(1,0){2}}
\put(10.2,4.1){\line(1,0){2}}\put(10.2,4.3){\line(1,0){2}}
\put(10.2,.2){\line(1,0){2}}\put(14.2,2.2){\line(1,0){2}}

\put(2.2,2.2){\line(1,1){2}}\put(2.2,2.2){\line(1,-1){2}}
\put(6.2,4.2){\line(1,1){2}}\put(6.2,4.2){\line(1,-1){2}}
\put(6.2,.1){\line(1,1){2}}\put(6.2,.1){\line(1,-1){2}}
\put(6.2,.3){\line(1,1){2}}\put(6.2,.3){\line(1,-1){2}}
\put(8.2,2.1){\line(1,-1){2}}\put(8.2,2.3){\line(1,-1){2}}
\put(8.2,6.2){\line(1,-1){2}}
\put(8.2,2.2){\line(1,1){2}}
\put(8.2,-1.7){\line(1,1){2}}\put(8.2,-1.9){\line(1,1){2}}
\put(12.2,4.2){\line(1,-1){2}}\put(12.2,.2){\line(1,1){2}}

\put(4.1,4.2){\line(1,-2){2}}\put(4.3,4.2){\line(1,-2){2}}
\put(4.1,.2){\line(1,2){2}}\put(4.3,.2){\line(1,2){2}}\put(4.2,.2){\line(1,2){2}}
\put(10.1,4.2){\line(1,-2){2}}\put(10.3,4.2){\line(1,-2){2}}\put(10.2,4.2){\line(1,-2){2}}
\put(10.1,.2){\line(1,2){2}}\put(10.3,.2){\line(1,2){2}}

\put(-.7,1.4){$182$}\put(1.2,1.4){$182$}\put(3.5,-.5){$82$}
\put(3.3,4.6){$100$}\put(5.5,-.5){$34$}\put(5.5,4.6){$16$}
\put(8,-2.6){$6$}\put(7,2){$11$}\put(8,6.5){$5$}
\put(10.1,-.6){$3$}\put(10.1,4.6){$5$}\put(12,-.6){$1$}
\put(12,4.6){$1$}\put(14,1.4){$1$}\put(16,1.4){$1$}
\end{picture}

\centerline{\scriptsize{6. The degrees of the Schubert classes}}

\bigskip

\medskip
For completeness we compile the other entries of the multiplication table. 
(The equivariant version could also be derived.)
\medskip

\begin{center}
\begin{tabular}{llllll}
$(\sigma_2)^2$ & = & $\sigma_4+2\sigma'_4+2\sigma''_4$\hspace*{2cm} & & & \\
$\sigma'_2\sigma_2$ & = & $\sigma_4+3\sigma'_4+\sigma''_4$ &
$(\sigma'_2)^2$ & = & $3\sigma_4+3\sigma'_4+\sigma''_4$  \\
$\sigma_3\sigma_2$ & = & $3\sigma_5+\sigma'_5$ & 
$\sigma_3\sigma'_2$ & = & $5\sigma_5+\sigma'_5$ \\
$\sigma'_3\sigma_2$ & = & $\sigma_5+\sigma'_5$ & 
$\sigma'_3\sigma'_2$ & = & $\sigma_5+\sigma'_5$ \\
$\sigma_4\sigma_2$ & = & $\sigma_6+\sigma'_6$ & 
$\sigma_4\sigma'_2$ & = & $\sigma_6+3\sigma'_6$ \\
$\sigma'_4\sigma_2$ & = & $2\sigma_6+3\sigma'_6$ & 
$\sigma'_4\sigma'_2$ & = & $3\sigma_6+3\sigma'_6$ \\
$\sigma''_4\sigma_2$ & = & $2\sigma_6+\sigma'_6$ & 
$\sigma''_4\sigma'_2$ & = & $\sigma_6+\sigma'_6$ \\
$\sigma_5\sigma_2$ & = & $\sigma_7$ & 
$\sigma_5\sigma'_2$ & = & $2\sigma_7$ \\
$\sigma_5\sigma_2 $& = & $\sigma_7$& 
$\sigma_5\sigma'_2 $& = & $2\sigma_7 $\\
$\sigma_6\sigma_2$ & = & $\sigma_8$ & 
$\sigma_6\sigma'_2$ & = & $0$ \\
$\sigma'_6\sigma_2$ & = & $0$ & 
$\sigma'_6\sigma'_2$ & = & $\sigma_8$ \\
 &&&&& \\
$(\sigma_3)^2$ & = & $3\sigma_6+5\sigma'_6$\hspace*{2cm} & & & \\
$\sigma'_3\sigma_3$ & = & $\sigma_6+\sigma'_6$ &
$(\sigma'_3)^2$ & = & $\sigma_6+\sigma'_6$  \\
$\sigma_4\sigma_3$ & = & $2\sigma_7$ & 
$\sigma_4\sigma'_3$ & = & $0$ \\
$\sigma'_4\sigma_3 $& = & $2\sigma_7$& 
$\sigma'_4\sigma'_3 $& = & $\sigma_7 $\\
$\sigma''_4\sigma_3 $& = & $0$& 
$\sigma''_4\sigma'_3 $& = & $\sigma_7 $\\
$(\sigma_4)^2$ & = & $\sigma_8$ & $\sigma_4\sigma'_4$ & = & $0$ \\ 
$(\sigma_4)^2$ & = & $\sigma_8$ & $\sigma_4\sigma''_4$ & = & $0$ \\ 
$(\sigma''_4)^2$ & = & $\sigma_8$ & $\sigma'_4\sigma''_4$ & = & $0$ 
\end{tabular}
\end{center}

\medskip
We can deduce a presentation of the cohomology ring over $\QQ$, choosing the hyperplane class
and the codimension two class $\sigma_2$ as generators. 

\begin{proposition}
The rational cohomology ring of $CG$ is 
$$H^*(CG,\QQ)=\QQ[\sigma_1,\sigma_2]/\langle \sigma_1^5-5\sigma_1^3\sigma_2+6\sigma_1\sigma_2^2,
16\sigma_2^3-27\sigma_1^2\sigma_2^2+9\sigma_1^4\sigma_2\rangle.$$
\end{proposition}

\subsection{The restriction map}

For future use we need to analyse the restriction of cohomology classes from $G=G(4,7)$ to $CG$. 
This is rather straightforward and can be done in different ways. In most cases one can simply take a 
Schubert class on $G$ and compute the degree of its intersection with the fundamental class of $CG$. 
If the result can be expressed uniquely as a non negative combination of the degrees of the 
Schubert classes of $CG$ of the same codimension, the relation we get between the degrees is 
also a relation between the classes, and we are done. When some ambiguity remains, we can 
compute the intersection of a general Schubert variety on $G$ in the Schubert class of interest, 
with the Schubert classes on $CG$. Since we know Poincar\'e duality on $CG$, the result follows. 

\begin{proposition}
The restriction map $\iota^* : A^*(G)\rightarrow A^*(CG)$ is given as follows, where we denote
by $\tau_\lambda$ the Schubert class on $G$ defined by a partition $\lambda$:
$$\begin{array}{ccccccccc}
\iota^*\tau_{1} &=& \sigma_1 & \iota^*\tau_{3} &=& \sigma'_3 & \iota^*\tau_{31} &=& \sigma'_4+ \sigma''_4 \\
\iota^*\tau_{2} &=& \sigma'_2 & \iota^*\tau_{21} &=& \sigma_3+ 2\sigma'_3 & \iota^*\tau_{22} &=& \sigma_4
 + \sigma'_4+ \sigma''_4 \\
\iota^*\tau_{11} &=& \sigma_2 & \iota^*\tau_{111} &=& \sigma_3 & \iota^*\tau_{211} &=& \sigma_4+2\sigma'_4 \\
 & & &  &&  & \iota^*\tau_{1111} &=& \sigma_ 4\\
\iota^*\tau_{32} &=& \sigma_5+\sigma'_5 & \iota^*\tau_{33} &=& \sigma_6+\sigma'_6  &  &&  \\
\iota^*\tau_{311} &=& \sigma_5+\sigma'_5  & \iota^*\tau_{321} &=& 3\sigma_6+3\sigma'_6 & \iota^*\tau_{222} &=& 2\sigma_6+2\sigma'_6   \\
\iota^*\tau_{221} &=& 3\sigma_5+\sigma'_5  & \iota^*\tau_{3111} &=& \sigma_6+\sigma'_6  &  \iota^*\tau_{2211} &=& \sigma_6+3\sigma'_6  \\
\iota^*\tau_{2111} &=& 2\sigma_5 &  &  && &&\\
&& &\iota^*\tau_{331} &=& 2\sigma_7 & \iota^*\tau_{332} &=& \sigma_8   \\
&&&\iota^*\tau_{322} &=& 2\sigma_7 & \iota^*\tau_{3311} &=& \sigma_8  \\
&&&\iota^*\tau_{3211} &=& 2\sigma_7 & \iota^*\tau_{3221} &=& \sigma_8   \\
&&&\iota^*\tau_{2221} &=& 2\sigma_7 & \iota^*\tau_{2222} &=& \sigma_8   \\
\end{array}$$
\end{proposition}

\medskip

\begin{corollary}
The image of $\iota^*$ is a sublattice of $A^*(CG)$ of index $16$.
\end{corollary} 
 
\subsection{The projective dual of $CG$}

As a byproduct of our computations we can derive interesting information on the projective
dual of $CG$, which is the variety $CG^{\vee}$ parametrizing, in the dual projective 
space, the tangent hyperplanes to $CG$. We will use the Katz-Kleiman formula \cite{gkz}, 
following which we should consider the polynomial
$$c_{CG}(q)=\sum_{i=0}^8q^{i+1}\int_{CG}c_{8-i}(\Omega_{CG})\sigma_1^i.$$
Then, if $c'_{CG}(1)\ne 0$, the dual variety of  $CG$ is a hypersurface of precisely 
that degree. 

\begin{proposition}
The projective dual variety $CG^{\vee}$ is a hypersurface of degree $17$ in $\PP^{28}$.
\end{proposition}

{\it Proof.} 
Having computed the weights of the tangent spaces at the $T$-fixed points,
we immmediately deduce the equivariant total Chern class of the tangent 
bundle $T_{CG}$. A routine computation then yields
$$\begin{array}{lll}
c_1(T_{CG}) & = & 4\sigma_1, \\
c_2(T_{CG}) & = & 9\sigma_2+ 7\sigma'_2,\\
c_3(T_{CG}) & = & 28\sigma_3+ 52\sigma'_3,\\
c_4(T_{CG}) & = & 49\sigma_4+ 88\sigma'_4+46\sigma''_4, \\
c_5(T_{CG}) & = & 76\sigma_5+ 160\sigma'_5,\\
c_6(T_{CG}) & = & 133\sigma_6+ 151\sigma'_6,\\
c_7(T_{CG}) & = & 90\sigma_7,\\
c_8(T_{CG}) & = & 15\sigma_8.
\end{array}$$
Using our computations of the degrees of the Schubert classes we deduce that
$$c_{CG}(q)=15q-90q^2+344q^3-860q^4+1492q^5-1784q^6+1438q^7-738q^8+182q^9.$$
This gives  $c'_{CG}(1)=17$ and the result follows.
$\Box$ 

\subsection{Concluding remarks} The cohomology of $CG$ looks very much like that of the
Grassmannian $G(2,6)$, which has the same Betti numbers and a similarly looking 
presentation, with relations in degree $5$ and $6$. An important difference between
the two varieties is that the index of $CG$ is four, while the index of $G(2,6)$ is six. This makes much harder
the task of computing the quantum cohomology ring of $CG$, since we have much more 
freedom in the possible quantum deformations of the two relations than in the case
of $G(2,6)$ (recall that for any prime Fano manifold, the 
degree of the quantum parameter is equal to the index).
Moreover the fact
that $CG$ is only quasi-homogeneous, not homogeneous, makes the computation of the 
Gromov-Witten invariants much less straightforward since they are not clearly enumerative. We plan to tackle these issues in a sequel to the present paper.

\bigskip 

{\scriptsize
{\sc Institut de Mathématiques de Marseille,  UMR 7373 CNRS/Aix-Marseille Universit\'e, 
Technop\^ole Château-Gombert, 
39 rue Frédéric Joliot-Curie,
13453 MARSEILLE Cedex 13,
France}

{\it Email address}:  {\tt laurent.manivel@math.cnrs.fr}


\begin{thebibliography}{aa}

\bibitem[Baez]{baez} Baez  J.,  {\it The octonions}, Bull. AMS {\bf 39} (2002), no. 2, 145-205.

\bibitem[BB]{bb}Bialynicki-Birula A., {\it Some theorems on actions of algebraic groups}, Annals of Math. (2) {\bf 98} (1973), 480-497.

\bibitem[Br]{br} 
Brion  M., {\it 
Equivariant cohomology and equivariant intersection theory},
in  Representation theories and algebraic geometry (Montreal 1997), 1–37, Kluwer 1998. 

\bibitem[BI]{bi} Brion M., Inamdar S.P., {\it Frobenius splitting of spherical varieties}, in Algebraic groups and their generalizations: classical methods, 207-218, Proc. Sympos. Pure Math. {\bf 56}, AMS 1994. 

\bibitem[GKZ]{gkz}
Gel'fand I.M., Kapranov M.M., Zelevinsky A.V.,  Discriminants, resultants, and multidimensional determinants,  Birkh\"auser 1994.

\bibitem[Kr]{kr} Kr\"amer M., {\it 
Sph\"arische Untergruppen in kompakten zusammenh\"angenden Liegruppen}, 
Compositio Math. {\bf 38} (1979), no. 2, 129–153. 

\bibitem[LM1]{lm1} Landsberg J.M., Manivel L., {\it 
The projective geometry of Freudenthal's magic square},
J. Algebra {\bf 239} (2001), no. 2, 477–512. 

\bibitem[LiE]{lie} LiE, A computer algebra package for Lie group computations, available at http://young.sp2mi.univ-poitiers.fr/~marc/LiE/. 


\end{thebibliography}
\end{document}